\documentclass[11pt]{article}

\usepackage[T2A]{fontenc}
\usepackage{amsmath}
\usepackage{amssymb,latexsym}
\usepackage{oldlfont}
\usepackage{amsfonts}

\usepackage[russian,english]{babel}

\newtheorem{theorem}{Theorem}[section]
\newtheorem{lemma}{Lemma}[section]
\newtheorem{pr}{Proposition}[section]
\newtheorem{cor}{Corollary}[section]
\newtheorem{df}{Definition}[section]

\newtheorem{ex}{Example}[section]

\renewcommand{\theequation}{\arabic{section}.\arabic{equation}}

\begin{document}

\renewcommand{\abstractname}{}

\renewcommand{\thesection}{\arabic{section}.}
\renewcommand{\theequation}{\arabic{section}.\arabic{equation}}



\vspace{5mm}

\begin{center}

{\Large \bf GEOMETRICAL AND ANALYTICAL CHARACTERISTIC PROPERTIES OF
PIECEWISE AFFINE MAPPINGS}

\end{center}

\begin{center}
{\bf V.V. Gorokhovik}
\end{center}

\begin{center}\textit{Institute of Mathematics \\ The National Academy of Sciences of Belarus \\
Surganova st., 11, Minsk, 220072, Belarus \\ e-mail:
gorokh@im.bas-net.by}\end{center}

 \setcounter{page}{1}

\vspace{7mm}

{\small {{\bf Abstract}
\smallskip

Let $X$ and $Y$ be finite dimensional normed spaces, ${\mathcal{F}}
(X,Y)$ a collection of all mappings from $X$ into $Y.$ A mapping $P
\in {\mathcal{F}} (X,Y)$ is said to be piecewise affine if there
exists a finite family of convex polyhedral subsets covering $X$ and
such that the restriction of $P$ on each subset of this family is an
affine mapping. In the paper we prove a number of characterizations
of piecewise affine mappings. In particular we show that a mapping
$P:X \to Y$ is piecewise affine if and only if for any partial order
$\preceq$ defined on $Y$ by a polyhedral convex cone both the
$\preceq$-epigraph and the $\preceq$-hypograph of $P$ can be
represented as the union of finitely many convex polyhedral subsets
of $X \times Y$. When the space $Y$ is ordered by a minihedral cone
or equivalently when $Y$ is a vector lattice the collection
${\mathcal{F}} (X,Y)$ endowed with standard pointwise algebraic
operations and the pointwise ordering is a vector lattice too. In
the paper we show that the collection of piecewise affine mappings
coincides with the smallest vector sublattice of ${\mathcal{F}}
(X,Y)$ containing all affine mappings. Moreover we prove that each
convex (with respect to an ordering of $Y$ by a minihedral cone)
piecewise affine mapping is the least upper bound of finitely many
affine mappings. The collection of all convex piecewise affine
mappings is a convex cone in ${\mathcal{F}} (X,Y)$ the linear
envelope of which coincides with the vector subspace of all
piecewise affine mappings.
 }

\vspace{3mm}


\vspace{3mm}


}
%
%

\section{Introduction}

Roughly speaking, a mapping is piecewise affine if it is ``glued''
with finitely many ``pieces'' of affine mappings (for an exact
definition see Section~3). By simplicity, this class of nonlinear
mappings is the most close one to linear and affine mappings.
Along with this, every continuous nonlinear mapping can be
approximated on any compact set by piecewise affine mappings with
an arbitrary accuracy. Due to such good approximating properties
piecewise affine mappings are widely used in nonlinear analysis
both in pure theoretical studies and in various applications (see,
for instance, [1~--~9]). In the present paper we review the
results of a number of articles [10~--~15] written by the author
in the recent past and devoted to piecewise affine functions and
mappings acting between finite-dimensional vector spaces. These
articles were published mainly in not readily available issues in
Russian. The main purpose of the paper is to make these results
more available to the interested readers. In Section~2 we present
some preliminary facts on polyhedral (not necessarily convex)
sets. Various definitions of piecewise affine mappings and some of
their properties are discussed in Section~3. Geometrical
characteristic properties of piecewise affine mappings are
presented in Section 4. In particular, it is proved that a mapping
$P:X \to Y$ is piecewise affine if and only if for any partial
order $\preceq$ defined on $Y$ by a polyhedral convex cone both
the $\preceq$-epigraph and the $\preceq$-hypograph of $P$ can be
represented as the union of finitely many convex polyhedral
subsets of $X \times Y$. In Section 6 for the case when $Y$ is a
vector lattice a number of analytical representations of piecewise
affine mappings are presented. The main above results are
specified in Section 7 to piecewise linear mappings.


\section{Preliminaries on polyhedral (not necessarily
convex) sets}

\setcounter{equation}{0}

Let $X$ be a finite-dimensional normed space and let $X^*$ be its
dual space whose elements are linear functions on $X.$

A hyperplane in $X$ is the set $H(a^*,\alpha ): = \{ x \in X\mid
a^*(x) = \alpha  \} ,$ where $a^*\in X^*,$ $a^*\ne 0,$ $\alpha
\in{\mathbb{R}}.$ Each hyperplane generates in $X$ two closed
halfspaces
$$H_\le(a^*,\alpha ): =
\{ x \in X\mid a^*(x) \le \alpha  \}\,\,
\text{and}\,\,H_\ge(a^*,\alpha ): = \{ x \in X\mid a^*(x) \ge \alpha
\}
$$
and two complementary open halfspaces
$$H_>(a^*,\alpha ): =
\{ x \in X\mid a^*(x) > \alpha  \} \,\,\text{and}\,\, H_<(a^*,\alpha
): = \{ x \in X\mid a^*(x) < \alpha  \}.$$

Since $H_\le(a^*,\alpha ) = H_\ge(-a^*,-\alpha )$ and
$H_>(a^*,\alpha )= H_<(-a^*,-\alpha ),$ we shall mainly present
closed halfspaces in the form $H_\le(a^*,\alpha )$ and open
halfspaces in the form $H_>(a^*,\alpha ).$

\vspace{1mm}

A convex subset $Q$ of $X$ is said to be {\it polyhedral}
[16~--~19] if it is the intersection of finitely many closed
halfspaces.  In other words a convex set $Q  \subset X$ is
polyhedral if it can be presented in the form $Q =
\bigcap\limits_{j = 1}^k {H_ \le (a_j^*,\alpha _j )},$ where
$a_j^*\in X^*,$ $\alpha _j \in{\mathbb{R}},$ $j=1,2,\ldots,k.$ By
the convention the whole space $X$ is a polyhedral convex set too.

\vspace{1mm}

Since every convex polyhedral set is in fact a solution set of
some finite system of linear inequalities, the theory of convex
polyhedral sets is substantially developed as a part of the theory
of linear inequalities [20, 21] and as a part of the linear
programming theory [22,23] which studies problems of minimizing
linear functions on convex polyhedral sets.

\vspace{1mm}

An arbitrary (not necessarily convex) set $Q \subset X$ is called
{\it polyhedral} [24] if it is the union of finitely many convex
polyhedral sets. (In [7] such sets were called piecewise
polyhedral.)

\vspace{1mm} It follows immediately from the above definition that
any polyhedral set $Q$ can be represented in the form
\begin{equation}\label{e2.1}
Q = \bigcup\limits_{i = 1}^m \bigcap\limits_{j = 1}^{k(i)} H_ \le
(b^* _{ij} ,\beta _{ij} ) ,
\end{equation}
where $b^* _{ij}\in X^*,$ $\beta _{ij} \in{\mathbb{R}},$
$j=1,\ldots,k(i);$ $i=1,\ldots,m.$

\vspace{2mm}

Let ${\mathcal M}(X)$ be the Boolean lattice of all subsets of $X$
with set operations of union and intersection as lattice
operations and let $M(X)$ be the collection of all polyhedral
subsets of $X.$ It was shown in [11, 12] that $M(X)$ is the
smallest sublattice in ${\mathcal M}(X),$ which contains all
closed halfspaces of the space $X.$  Thus polyhedral sets of $X$
are exactly those which can be constructed from a finite family of
closed halfspaces as a result of finitely many operations of union
and intersection. It is follows from characteristics of general
sublattice  (see [25, 26]), that every polyhedral set $Q$ can be
represented in the alternative form
\begin{equation}\label{e2.2}
Q = \bigcap\limits_{j = 1}^{k}\bigcup\limits_{i = 1}^{m(j)}  H_
\le (c^* _{ij} ,\gamma _{ij} ) ,
\end{equation}
where $c^* _{ij}\in X^*,$ $\gamma _{ij} \in{\mathbb{R}},$
$i=1,\ldots,m(j);$ $j=1,\ldots,k.$

In general the collections of closed halfspaces $H_ \le (b^* _{ij}
,\beta _{ij} ) , j=1,\ldots,k(i),\,i=1,\ldots,m,$ and $H_ \le (c^*
_{ij} ,\gamma _{ij} ) , i=1,\ldots,m(i),\,j=1,\ldots,k,$ which are
used for the representation of the same set $Q$ in \eqref{e2.1} and
\eqref{e2.2} can be different.

\vspace{1mm}

It was shown in [14] that for every polyhedral set $Q \subset X$
there exists a finite two-index family of closed halfspaces $\{ H_
\le (a^* _{ij} ,\alpha _{ij} ),$ $ i=1,\ldots,m; j=1,\ldots,k\}$
such that
$$
P = \bigcup\limits_{i = 1}^m \bigcap\limits_{j = 1}^{k} H_ \le (a^*
_{ij} ,\alpha _{ij} ) = \bigcap\limits_{j = 1}^{k} \bigcup\limits_{i
= 1}^m H_ \le (a^* _{ij} ,\alpha _{ij} ). $$

\vspace{2mm}


\section{The definition and elementary properties of
piecewise \\ affine mappings}

\setcounter{equation}{0}

Let $X$ and $Y$ be finite-dimensional normed spaces..

A finite family $\sigma = \{M_1, \ldots, M_k\}$ of convex polyhedral
subsets $M_1, \ldots, M_k$ of $X$ is called \textit{a polyhedral
covering} of a polyhedral set $Q \subset X,$ if
\begin{equation}
  \label{e1.1}
M_i \subset Q,\, i=1,\ldots,k; \,\,\text{and}\,\,
Q=\bigcup\limits_{i=1}^kM_i.
\end{equation}

A family $\sigma = \{M_1, \ldots,  M_k\}$ is called {\it a
polyhedral partition} of a polyhedral set $Q \subset X,$ if it is a
polyhedral covering and, in addition, the conditions
${\rm ri} M_i\cap {\rm ri} M_j=\varnothing,\ i, j = 1, \ldots, k,\
i\ne j,$
hold.

Here ${\rm ri} M$ stands for the relative interior of a convex set
$M.$

\vspace{3mm}

A polyhedral covering (a polyhedral partition) $\sigma = \{M_1,
\ldots,  M_k\}$ of the set $Q$ is called {\it solid}, if each $M_i,\
i=1, \ldots, k$ has nonempty interior or, equivalently, if
$$\dim{\rm aff} M_i=\dim X,\ i = 1, \ldots, k,$$
where ${\rm aff} M$ is the affine hull of the set $M.$

\begin{df}{\rm [13, 15].}\label{df1}
{\rm A mapping $P:X\to Y$ is said to be {\it piecewise affine}, if
there exists a polyhedral covering $\sigma = \{M_1, \ldots, M_k\}$
of the vector space $X$ and the collection of affine mappings
\linebreak $A_i:X\to Y, i=1, \ldots, k,$ such that
$$
P(x)=A_i(x),\  x\in M_i,\  i=1, \ldots, k.
$$}
\end{df}

In what follows we shall denote the collection of all piecewise
affine mappings from $X$ into $Y$ by the symbol $PA(X, Y).$ When
$Y={\mathbb{R}}$ the collection $PA(X, {\mathbb{R}})$ of all
piecewise affine functions from $X$ into ${\mathbb{R}}$ will be
denoted, for short, by $PA(X).$

\begin{ex}{\bf \hspace{-6pt}.}\label{ex1}
{\rm Every affine and, consequently, linear mapping is piecewise
affine. Hence, the following  inclusions hold:
$$
L(X, Y)\subset A(X, Y)\subset PA(X, Y).
$$
where $L(X, Y)$ and $A(X, Y)$ are, respectively, the space of
linear mappings and the space of affine ones from $X$ into $Y.$}
\end{ex}

\begin{ex}{\bf \hspace{-6pt}.}\label{ex2} {\rm The mappings
$F:{\mathbb{R}}\times{\mathbb{R}}\to{\mathbb{R}}$ and
$G:{\mathbb{R}}\times{\mathbb{R}}\to{\mathbb{R}}$, defined by
$F(x_1, x_2)=\max\{x_1, x_2\}$ and $G(x_1, x_2)=\min\{x_1, x_2\},$
are piecewise affine. The polyhedral covering of $X$ which is
associated with $F$ and $G$ is $\{M_1, M_2\}$, where $M_1:=\{(x_1,
x_2)\in{\mathbb{R}}\times{\mathbb{R}}\mid x_1\le x_2\}$,
$M_2:=\{(x_1, x_2)\in{\mathbb{R}}\times{\mathbb{R}}\mid x_1\ge
x_2\}$.}
\end{ex}

\begin{ex}{\bf \hspace{-6pt}.}\label{ex3} {\rm The mapping $G_\lambda: X \ni x \to
\lambda x \in X$, where $\lambda$ is an arbitrary fixed real number,
and the mapping $F: X\times X \ni (x_1, x_2)\to x_1+x_2 \in X$ are
linear and, consequently, piecewise affine.}
\end{ex}

The above definition of piecewise affine mappings differs from the
commonly accepted one (see, for instance, [4, 5]) in the
requirement that a family $\sigma$ is a covering of $X$ rather
than a solid partition of $X.$ This requirement is less
restrictive and therefore it is more convenient in proving some
properties of piecewise affine mappings. Moreover, the above
definition is actually equivalent to the common one. For the first
time, this fact was proved in [13]. Later, another proof was given
in [7]. Below, we will reproduce the proof given in [13].

\begin{pr}{\bf \hspace{-6pt}.}\label{pr3.1}
If $\sigma = \{M_1, \ldots, M_k\}$ is a polyhedral covering of a
convex polyhedral set $Q$ with ${\rm int}Q \ne \varnothing$ then
the family
$$
\sigma'=\{M_i\in\sigma\mid {\rm int} M_i \ne \varnothing\}
$$
is a solid polyhedral covering of the set $Q.$
\end{pr}

\noindent{\bf Proof.} First we prove the proposition for the case
when $Q = X.$  Let ${\hat \sigma}:= \{M_i\in\sigma\mid {\rm int}
M_i =\varnothing\}$. The set $\hat M:=\cup\{M\mid
M\in\widehat\sigma\}$ is closed and nowhere dense in $X$ (as the
union of finitely many closed and nowhere dense sets [27, p. 114,
Proposition 1]), and, hence, its complement $X\backslash {\hat M}$
is everywhere dense in $X.$ Then, since $X\backslash{\hat
M}\subset M',$ the set $M':=\cup\{M\mid M\in\sigma'\}$ is dense in
$X$ too. Because $M'$ is closed we get $M'=X$. This proves the
assertion of the proposition for the case $Q=X.$

Now we suppose that $\sigma=\{M_1, \ldots, M_k\}$ is a polyhedral
covering of an arbitrary convex polyhedral set $Q$ with ${\rm
int}Q$ and $\sigma'=\{M\in\sigma\mid{\rm int} M\ne\varnothing\}.$
We can present the set $Q$ in the form
$Q=\bigcap\limits_{i=1}^pH_{\le}(a^*_i, \alpha_i)$, where
$a^*_i\in X^*\backslash\{0\},\ \alpha_i\in{\mathbb{R}},\
i=1,\ldots,p.$. It is not difficult to see, that the family
$\sigma\cup\{H_{\ge}(a^*_1, \alpha_1), \ldots, H_{\ge}(a^*_p\,,
\alpha_p)\}$ is a polyhedral covering of the whole space $X.$
Then, as it has been proved in the first part of this proof, the
family $\sigma'\cup\{H_{\ge}(a^*_1, \alpha_1), \ldots,
H_{\ge}(a^*_p\,, \alpha_p)\}$ is a solid polyhedral covering of
$X$. Now we conclude from the equality ${\rm int} Q =
\bigcap\limits_{i=1}^pH_<(a^*_i, \alpha_i)$ (see [16, Theorem
6.5]) that ${\rm int} Q \subset \cup\{M\mid M\in\sigma'\} \subset
Q$ and, since the set $Q$ is convex and closed, we have ${\rm cl}
({\rm int} Q)= Q$ (see., for instance, [16, Theorem 6.3]) and,
consequently,
 $Q \subset {\rm cl}(\cup\{M\mid M\in \sigma'\})=\cup\{M\mid M\in\sigma'\} \subset Q$.

It is easy to see that when a family $\sigma$ is a polyhedral
partition of a set $Q$, then $\sigma'$ is a solid polyhedral
partition of $Q.$

This completes the proof of Proposition \ref{pr3.1}.

\begin{pr}{\bf \hspace{-6pt}.}\label{pr3.2}
For any polyhedral covering $\sigma = \{M_1, \ldots, M_k\}$ of a
polyhedral set $Q$ there exists a polyhedral partition $\omega =
\{D_1, \ldots, D_m\}$ of the set $Q$ such that every $D_j \in
\omega$ is contained in some $M_i\in \sigma.$
\end{pr}

First of all we prove the following auxiliary lemma.

\begin{lemma}{\bf \hspace{-6pt}.}\label{l3.1}
Let $\{H(a^*_i, \alpha_i),\ i\in S\}$ be a family of hyperplanes
in $X$ indexed by elements of a set $S$ and let the set
$$
D_I=(\bigcap_{i\in I}H_{\le}(a_i^*, \alpha_i))\bigcap
    (\bigcap_{i\in S\backslash I}H_{\ge}(a^*_i, \alpha_i))
$$
be associated with every (possibly, empty) subset $I$ of $S$.

Then for any subsets $I,\,J \subset S$ one has either $D_I = D_J$ or
${\rm ri}D_I \cap {\rm ri}D_J = \varnothing.$
\end{lemma}

\noindent{\bf Proof of Lemma \ref{l3.1}.} Notice that for every
$i\in K:=(I\backslash J)\bigcup(J\backslash I)$ the sets $D_I$ and
$D_J$ lies in different closed halfspaces generated by the
hyperplane $H(a^*_i, \alpha_i).$

When there is $i\in K$ such that at least one of the sets $D_I$ or
$D_J$ does not lie whole in $H(a^*_i, \alpha_i),$ the convex sets
$D_I$ and $D_J$ are properly separated by the hyperplane $H(a^*_i,
\alpha_i)$ and, hence, ${\rm ri} D_I\bigcap {\rm ri}
D_J=\varnothing$ (see. [16, Theorem 11.3]).

In the case when $D_I,\,D_J \subset H(a^*_i, \alpha_i),$ for all
$i\in K$, we have
$$
D_I=(\bigcap_{i\in K}H(a_i^*, \alpha_i))\bigcap
    (\bigcap_{i\in I\backslash K}H_{\le}(a^*_i, \alpha_i))\bigcup
    (\bigcup\limits_{i\in S\backslash(I\cup K)}H_{\ge}(a^*_i, \alpha_i)),
$$
$$
D_J=(\bigcap_{i\in K}H(a_i^*, \alpha_i))\bigcap
    (\bigcap_{i\in J\backslash K}H_{\le}(a^*_i, \alpha_i))\bigcup
    (\bigcup\limits_{i\in S\backslash(J\cup K)}H_{\ge}(a^*_i, \alpha_i)).
$$
Since $I\backslash K=J\backslash K= I\cap J$ and $I\cup K= J\cup K=
I\cup J$, then $D_I=D_J$. This proves the lemma.

\medskip\noindent{\bf Proof of Proposition \ref{pr3.2}.} Let $\sigma = \{M_1,
\ldots, M_k\}$ be a polyhedral covering of the set $Q.$ Every set
$M_i, i=1, \ldots, k$ can be presented in the form
$M_i=\bigcap\limits_{j\in S(i)}H_{\le}(a^*_j, \alpha_j)$, where
$S(i)$ is a finite family of indices. Let $\{H_{\le}(a^*_j,
\alpha_j),\,j\in S:=\bigcup\limits_{i=1}^mS(i)\}$ be the collection
of all hyperplanes which are used in the above presentations of the
sets $\{M_1, \ldots, M_k\}$. With every nonempty subset $I\subset S$
we associate the set $D_I=(\bigcap\limits_{j\in I}H_{\le}(a^*_j,
\alpha_j))\bigcap$ $\bigcap(\bigcap\limits_{j\in S\backslash
I}H_{\ge}(a^*_j, \alpha_j)).$ Notice that for some $I,\,J \subset
S,\, I \ne J,$ we can have $D_I =D_J.$  Let us consider the family
$\omega\,'$ consisting of all nonempty subsets $D_I$, that lie in
one of subsets  $M_i,\ i=1, \ldots, k$. It is evident that
$\cup\{D\mid D\in\omega\,'\}\subset Q$.

Let $x$ be a point of $Q$ and let $i\in\{1, 2, \ldots, m\}$ be
such that $x\in M_i.$ It is easy to see that $x\in D_J,$ where
$J=\{j\in S\mid x\in H_{\le}(a^*_j, \alpha_j)\}$. Moreover, since
$S(i)\subset J$, then $D_J\subset M_i$ and, hence, $D_J\in
\omega\,'$. As $x$ is an arbitrary point of $Q,$ we get
$Q\subset\bigcup\{D\mid D\in \omega\,'\}$. Consequently, the
family $\omega\,'$ is a covering of $Q.$ Besides it follows from
Lemma~1 that different subsets of $\omega\,'$ do not intersect
each other by their relative interiors. Thus, the subfamily
$\omega$ of nonempty different subsets of $\omega\,'$ is in fact a
polyhedral partition of $Q$ and, moreover, each subset of $\omega$
is contained in some subset of the covering $\sigma.$ This
completes the proof of Proposition \ref{pr3.2}.

\medskip
The next theorem follows immediately from Proposition \ref{pr3.1}
and Proposition \ref{pr3.2}.\medskip

\begin{theorem}{\bf \hspace{-6pt}.}\label{th3.1}
A mapping $P:X \to Y$ is piecewise affine if and only if there
exists a solid polyhedral partition $\sigma'$ of the space $X$ such
that $P$ coincides with some affine mapping on each subset of
$\sigma'.$
\end{theorem}

Theorem \ref{th3.1} proves that Definition \ref{df1} of piecewise
affine mappings is in fact equivalent to the conventional one.

\medskip

We complete this section with describing some properties of the
space of piecewise affine functions.

\begin{theorem}{\bf \hspace{-6pt}.}\label{th3.2}
The composition of piecewise affine mappings is a piecewise affine
mapping.
\end{theorem}

\noindent{\bf Proof.} Let $X, Y$ and $Z$ be finite dimensional
normed spaces and let $P\in PA(X, Y)$ and $Q\in PA(Y, Z)$ be
piecewise affine mappings. Let $\sigma=\{M_1, \ldots, M_k\}$ and
$\delta=\{D_1, \ldots, D_p\}$ be polyhedral covering of $X$ and
$Y$ which correspond, respectively, to $P$ and $Q;$  $\{P_1,
\ldots, P_k\}$ and $\{Q_1, \ldots, Q_p\}$ the collections of
affine mappings from $A(X, Y)$ and from $A(Y, Z)$ such that
$P(x)=P_i(x),\ x\in M_i,\ i=1, \ldots, k;\ Q(y)=Q_j(y),\ y\in
D_j,\ j=1, \ldots, p$. The sets $C_{ij}:=M_i\cap P^{-1}_i(D_j),\
i=1, \ldots, k,\ j=1, \ldots, p, $ are convex and polyhedral (it
follows from properties of convex polyhedral sets [13 -- 16]) and
form a polyhedral covering of $X$. Besides, on each set $C_{ij}$
the composition $Q\circ P$ coincides with the affine mapping
$Q_j\circ P_i.$ This proves the theorem.

Theorem \ref{th3.2} implies a number of important corollaries.

\begin{cor}{\bf \hspace{-6pt}.}\label{cor3.1}
The collection $PA(X, Y)$ endowed with standard pointwise
operations of addition and multiplication by reals is a vector
space.
\end{cor}

The validity of this assertion immediately follows from Theorem
\ref{th3.2} and Example \ref{ex3}.

\vspace{3mm}

Let ${\rm dim} Y = n$ and let $\{e_1, \ldots, e_n\}$ be a vector
basis of $Y.$ As is well known any mapping $P:X\to Y$ can be
uniquely associated with its {\it coordinate functions} $p_i:X \to
{\mathbb{R}},$ $i=1,\ldots, n$, such that
$P(x)=p_1(x)e_1+\ldots+p_n(x)e_n$ for all $x \in X.$

\begin{cor}{\bf \hspace{-6pt}.}\label{cor3.2}
The mapping $P:X\to Y$ is piecewise affine if and only if its
coordinate functions $p_i:X\to {\mathbb{R}},\ i=1,\ldots, n$, are
piecewise affine.
\end{cor}

\noindent{\bf Proof.} Since the mappings $\phi_i:\lambda \to
\lambda e_i,$ $i=1,\ldots, n,$ from $\mathbb{R}$ into $Y$ are
linear, the sufficient part of Corollary \ref{cor3.2} follows
directly from Theorem \ref{th3.2} and Corollary \ref{cor3.1}.

To prove the necessary part we choose in $Y^*$ the basis $\{e^*_1,
\ldots, e^*_n\}$ which is dual to the basis $\{e_1, \ldots, e_n\}$
of $Y.$  It means that $\langle e^*_i,\, e_i \rangle =1$ for all
$i=1,\ldots,n$ and $\langle e^*_i,\, e_j \rangle =0$ for
$i,j=1,\ldots,n;\,i \ne j.$ Since $e^*_i,\ i=1,\ldots, n,$ are
linear functions on $Y,$ it follows from Theorem \ref{th3.2} that
the functions $p_i(x)=e^*_i(P(x)),\, i=1, \ldots, n,$ are
piecewise affine. \vskip4pt plus2pt

%
%
%
%
%
%



\section{Geometrical characteristic properties of piecewise
affine mappings}

\setcounter{equation}{0}

Directly from the definition of piecewise affine mappings we can see
that the graph of a piecewise affine mapping is a union of finitely
many convex polyhedral subsets, each of which is a part of the graph
of an affine mapping, and, consequently, the graph of a piecewise
affine mapping is a polyhedral (but, in general, nonconvex) set. In
spite of apparent evidence of this observation we provide it with a
proof.

\begin{theorem}{\bf \hspace{-6pt}.}\label{th4.1}
A mapping $P:X\to Y$ is piecewise affine if and only if its graph
${\rm graph}P:=\{(x, y)\in X\times Y\mid P(x)=y\}$ is a polyhedral
set in $X\times Y$.
\end{theorem}

\medskip\noindent{\bf Proof.} {\it Necessity.} Let $P:X\to Y$ be a
piecewise affine mapping and let $\sigma=\{M_1, \ldots, M_k\}$ and
${\mathcal{A}}= \{A_1,\ldots, A_k\}$ be the polyhedral partition of
the space $X$ and the collection of affine functions associated with
$P.$ Then the equality
$$
{\rm graph}P=\bigcup_{i=1}^k\left((M_i\times Y)\bigcap {\rm
graph}A_i\right),
$$
holds and we see that the graph of $P$ is a polyhedral set in
$X\times Y.$

{\it Sufficiency.} First we consider the case when $Y=\mathbb{R}.$
Let ${\rm graph}P=\displaystyle\bigcup_{i=1}^kG_i$, where
$\Sigma=\{G_1, \ldots, G_k\}$ is a family of convex polyhedral
sets in $X\times \mathbb{R}$. Then the family $\sigma:=\{M_1,
\ldots, M_k\}$ with $M_i:={\rm pr}_XG_i$ (${\rm pr}_XG_i$ stands
for the projection of $G_i$ on $X$) is a polyhedral covering of
$X.$ In view of Proposition~\ref{pr3.1} the subcollection
$\sigma'=\{M_i\in \sigma\mid {\rm int} M_i\ne\varnothing\}$ is a
solid polyhedral covering of $X.$ Let $M_i$ be a subset of
$\sigma'$ and $G_i$ be a subset of $\Sigma$ corresponding to
$M_i.$ Suppose that $\dim X=m.$ Since ${\rm int}
M_i\ne\varnothing$ we can choose in $M_i$ a collection $\{x_0,
\ldots, x_m\} \subset M_i$ of affinely independent points [16]. It
is not difficult to verify that the collection $\{(x_0, P(x_0)),
\ldots, (x_m, P(x_m))\}$ of points of $G_i$ is also affinely
independent and, consequently, $\dim{\rm aff}\,G_i\ge m$. But the
dimension of the affine hull of $G_i$ can not be equal to $m+1$
because it contradicts the single-valuedness of $P$. Therefore,
$\dim{\rm aff}\,G_i=m$ and, hence, the affine hull ${\rm
aff}\,G_i$ is a hyperplane in $X\times \mathbb{R}.$  Then there
are $a_i^*\in X^*$, and $\alpha_i, \beta_i\in \mathbb{R}$ such
that ${\rm aff}\,G_i=\{(x, \xi)\in X\times {\mathbb{R}} \mid
{a^*_i(x)}+ \alpha_i\xi=\beta_i\}.$  Note that the coefficient
$\alpha_i$ is not equal to zero, otherwise we had the inclusion
$M_i\subset\{x\in X\mid a_i^*(x) =\beta_i\}$ that contradicts the
condition ${\rm int} M_i\ne\varnothing$. Hence the hyperplane
${\rm aff}\,G_i$ is the graph of the affine function $h_i:x \to
({a^*_i(x)}-\beta_i)/\alpha_i.$ Thus the solid polyhedral covering
$\sigma'=\{M_i\in \sigma\mid {\rm int} M_i\ne\varnothing\}$ of the
space $X$ is associated with the collection of affine functions
$h_i:X\to\mathbb{R}$ such that $P(x)=h_i(x)$ for all $x\in M_i$.
This proves that the function $P:X\to\mathbb{R}$ is piecewise
affine.

Suppose now that  $\dim Y=n>1$ and choose in $Y$ a vector basis
$\{e_1, \ldots, e_n\}.$ Let \linebreak $p_i:X \to {\mathbb{R}},\
i=1, \ldots, n,$ be coordinate functions of the mapping $P$
corresponding to this basis. Then ${\rm graph}p_i=r_i({\rm graph}P),
\ i=1, \ldots, n$, where  the mapping $r_i:X\times Y\to
X\times\mathbb{R}$ is defined by the equality $r_i(x, y)=(x,
{\langle y, e_i^* \rangle}).$ Here $\{e_1^*, \ldots, e^*_n\}$ stand
for the basis in $Y^*$ which is dual to the basis $\{e_1, \ldots,
e_n\}$. Since the mappings $r_i,\ i=1, \ldots, n,$ are linear and
${\rm graph}P$ is polyhedral, ${\rm graph}p_{\hspace{1pt}i}$ is
polyhedral too. Consequently, as it was proved above, the coordinate
functions $p_i,\ i=1, \ldots, n,$ are piecewise affine. Then it
follows from Corollary~\ref{cor3.2} that the mapping $P$ is
piecewise affine too. This completes the proof of the theorem.
\vskip4pt plus2pt

A partial order $\preceq$ defined on $Y$ will be called {\it
polyhedral} if it is compatible with algebraic operations on $Y$ and
its positive cone $Y^+_\preceq:=\{y\in Y\mid 0\preceq y\}$ is
polyhedral.

\begin{theorem}{\bf \hspace{-6pt}.}\label{th4.2}
A mapping $P:X\to Y$ is piecewise affine if and only if for any
polyhedral partial order $\preceq$ defined on $Y$ both the
$\preceq$-epigraph ${\rm epi}_\preceq P:=\{(x, y)\in X \times
Y\mid P(x)\preceq y \}$ and the $\preceq$-hypograph ${\rm
hyp}_\preceq P:=\{(x, y)\in X\times Y\mid y \preceq P(x) \}$ of
$P$ are polyhedral sets in $X\times Y.$
\end{theorem}

\noindent{\bf Proof.}  Since for any polyhedral partial order
$\preceq$ defined on $Y$ both $\preceq$-epigraph ${\rm
epi}_\preceq A$ and $\preceq$-hypograph ${\rm hyp}_\preceq A$ of
any affine mapping  $A:X \to Y$ are convex polyhedral sets in $X
\times Y$ (it follows, for instance, from [16, Theorem 3.6]), the
necessary part of the theorem follows from the equalities
$$
{\rm epi}_\preceq P = \{\bigcup\limits_{i=1}^k((M_i\times
Y)\bigcap{\rm epi}_\preceq A_i)\},
$$
$$
{\rm hyp}_\preceq P = \{\bigcup\limits_{i=1}^k((M_i\times
Y)\bigcap{\rm hyp}_\preceq A_i)\},
$$
where  $\{M_1, \ldots, M_k\}$ is a polyhedral covering of $X$, and
$\{A_1, \ldots, A_{\hspace{1pt}k}\}$ is a collection of affine
functions such that $p(x)=A_{\hspace{1pt}i}(x)$ for all  $x\in M_i,\
i=1, \ldots, k$.

The sufficient part follows from the equality ${\rm graph}P={\rm
epi}_\preceq P\cap{\rm hyp}_\preceq P$ and Theorem~\ref{th4.1}.

\begin{cor}{\rm [10]}.\label{cor4.1}
A real-valued function $p:X\to \mathbb{R}$ is piecewise affine if
and only if its epigraph ${\rm epi}p:=\{(x, \alpha)\in X \times
{\mathbb{R}}\mid p(x)\le\alpha\}$ and its hypograph ${\rm
hyp}p:=\{(x, \alpha)\in X\times{\mathbb{R}}\mid p(x)\ge\alpha\}$ are
polyhedral sets.
\end{cor}


\section{Analytical characterizations of piecewise
affine mappings}

\setcounter{equation}{0}

In what follows we shall suppose that the vector space $Y$ is
endowed with a partial order $\preceq$ with respect to which $Y$
is an Archimedean vector lattice [25,\,26]. As it is known (see,
for instance, [28]), a finite-dimensional ordered vector space
$(Y,\,\preceq)$ is an Archimedean vector lattice if and only if
its positive cone $Y^+_\preceq:=\{y\in Y\mid 0\preceq y\}$ is a
conic convex hull of some vector basis $\{e_1,\,e_2,\ldots,e_n\}$
of $Y,$ where $n= {\rm dim}Y.$ Thus, the assumption that $Y$ is an
Archimedean vector lattice does not restrict the generality of
consideration because it is equivalent to the choice of some
vector basis in $Y.$  Below without additional mentions we shall
assume that a vector bases $\{e_1, \ldots, e_n\}$ of $Y$ and a
partial order $\preceq$ defined on $Y$ are compatible in such way
that $Y^+_\preceq:= {\rm coconv}\{e_1, \ldots, e_n\}$, where ${\rm
coconv}M$ denotes a conic convex hull of a set $M.$

Under above assumptions every mapping $F:X\to Y$ can be represented
as $F(x)= f_1(x)e_1+f_2(x)e_2+\ldots+f_n(x)e_n,$ where $f_i:X \to
Y,\,i=1,\,2,\ldots,n,$ are real-valued functions called the
coordinate functions of $F.$

The collection ${\mathcal{F}}(X,\,Y)$ of all mappings from $X$ into
$Y,$ endowed with standard pointwise algebraic operations and with
the partial order $\preccurlyeq$ defined by
\begin{equation}
  \label{e5.1}
F\preccurlyeq G \Longleftrightarrow F(x)\preceq G(x)\,\, \forall x
\in X \Longleftrightarrow f_i(x)\le g_i(x)\,\,\forall x \in
X,\,i=1,\ldots,n,
\end{equation}
(here $\{f_1,\ldots,f_n\}$ and $\{g_1,\ldots,g_n\}$ stand for
coordinate functions of $F$ and $G,$ respectively)\\ is a vector
lattice with
$$\sup\{F, G\}(x)=
(\max\{f_1(x),\, g_1(x)\},\ldots,\max\{f_n(x),\, g_n(x)\})$$ and
$$\inf(F, G)(x)=(\min\{f_1(x),\, g_1(x)\},\ldots,\min\{f_n(x),\,
g_n(x)\})$$ as lattice operations.

It is not difficult to see that, whenever $F$ and $G$ are piecewise
affine mappings, $\sup\{F, G\}$ and $\inf(F, G)$ are also piecewise
affine ones. Thus the following assertion is true.

\begin{pr}{\bf \hspace{-6pt}.}\label{pr5.1}
The collection of all piecewise affine mappings $PA(X, Y)$ from
$X$ into $Y$ with pointwise algebraic operations and with the
partial order $\preccurlyeq$ defined by \eqref{e5.1} is a vector
sublattice of the vector lattice ${\mathcal{F}}(X,\,Y).$
\end{pr}

A mapping $F:X\to Y$ is said to be $\preceq-${\it convex} (convex
with respect to a partial order$\preceq$) if
$$
F(\alpha x+\beta y)\preceq \alpha F(x)+\beta F(y),
$$
for all $x, y\in X$ and all $\alpha, \beta \ge 0,\ \alpha+\beta=1$.

\begin{pr}{\bf \hspace{-6pt}.}\label{pr5.2}
A piecewise affine mapping $P:X\to Y$ is $\preceq-${\it convex} if
and only if it can be represented in the form
\begin{equation}
  \label{e5.2}
P(x)=\sup\limits_{j\in J}A_j(x), \ x\in X,
\end{equation}
where $\{A_j:X\to Y,\ j\in J\}$ is a finite family of affine
mappings from $X$ into $Y$ and $\sup$ is the least upper bound in
the vector lattice $Y.$
\end{pr}

\noindent{\bf Proof.}  Since $P$ is $\preceq$-convex and piecewise
affine, the coordinate functions $p_i:X \to Y,\,i=1,\ldots,n,$ of
$P$ are convex and piecewise affine. Due to Proposition 3.1 of [7]
every coordinate function \linebreak $p_{\hspace{1pt}i}:X\to Y \,
(i=1, \ldots, n),$ can be represented in the form
$p_{\hspace{1pt}i}(x)=\max\limits_{j\in J_i}(a^*_j(x)+\alpha_j),\,
x \in X,$ $i=1, \ldots, k,$ where $J_i$ are finite family of
indices, $a^*_j:X\to {\mathbb R},$ $j\in J_1\cup J_2\cup \ldots
\cup J_n,$ are real-valued linear functions. Define the set of
multi-indices $J=J_1\times J_2\times \ldots \times J_n$ and
associate with every$j=(j_1, \ldots, j_n)\in J$ the affine mapping
$A_j=(a^*_{j_{\hspace{1pt}1}}(\cdot)+\alpha_{j_{\hspace{1pt}1}},
\ldots,
a^*_{j_{\hspace{1pt}n}}(\cdot)+\alpha_{j_{\hspace{1pt}n}}):X \to
Y.$ It is not difficult to verify that $P(x)=\sup\limits_{j\in
J}A_j(x), \ x\in X,$ where $\sup$ is the least upper bound in the
vector \linebreak lattice $Y.$

The converse assertion follows from the inequality
$$
\sup\limits_{i\in J}(a_i+b_i) \le \sup_{i\in J}a_i+\sup_{i\in J}b_i
$$
where $\{a_i\mid i\in J\}$ and $\{b_i\mid i\in J\}$ are finite
family of vectors in $Y.$

This completes the proof.

\begin{theorem}{\bf \hspace{-6pt}.}\label{th5.1}
Let $P:X \to Y$ be a mapping from $X$ into $Y$.

The following assertions are equivalent:

a) $P:X\to Y$ is piecewise affine;

b) $P:X\to Y$ can be represented in the form
\begin{equation}
  \label{e5.3}
P(x) = \inf\limits_{1\le i\le k}\sup\limits_{1\le j\le
m(i)}A_{ij}(x),\ \ x\in X,
\end{equation}
where $A_{ij}:X\to Y,\ j=1, \ldots, m(i),\ i=1, \ldots, k$, are
piecewise affine mappings;

c) $P:X\to Y$ can be represented in the form
\begin{equation}
  \label{e5.4}
P(x) = \sup\limits_{1\le i\le k}\inf\limits_{1\le j\le
m(i)}B_{ij}(x),\ \ x\in X,
\end{equation}
where $B_{ij}:X\to Y,\ j=1, \ldots, m(i),\ i=1, \ldots, k$, are
piecewise affine mappings;


d) $P:X\to Y$ can be represented in the form
\begin{equation}
  \label{e5.5}
P(x)=\sup\limits_{1\le i \le k}C_i(x) - \sup\limits_{1\le j \le
m}D_j(x),\, x\in X,
\end{equation}
where $C_i, D_j:X\to Y,\ i=1,\ldots, k, \, j=1, \ldots, m,$ are
piecewise affine mappings.
\end{theorem}

\noindent{\bf Proof.} The implications $b)\Rightarrow a),$
$c)\Rightarrow a)$ and $d)\Rightarrow a)$ immediately follow from
the fact that the collection $PA(X,\,Y)$ of all piecewise affine
mappings is a vector lattice (see Proposition \ref{pr5.1}) and
that affine mappings belong to $PA(X,\,Y).$


Let $P:X\to Y$ be a piecewise affine mapping. It follows from
Corollary~\ref{cor3.2} that coordinate functions
$p_s:X~\to~{\mathbb{R}},$ $s=1, \ldots, n$, corresponding to $P$
also are peacewise affine. Due to Theorem 3.1 and Proposition 3.1
of [10] we can represent each function $p_s:X~\to~{\mathbb{R}},$
$s=1,\ldots,n$ in the form
$$
p_s(x)=\min\limits_{i\in I_s}q_i(x), \ x\in X,
$$
where $I_s$ is a finite family of indices and  $q_i:X\to {\mathbb
R},\ i\in I_s,$ are convex piecewise affine functions. For every
$i=(i_1, \ldots, i_n)\in I:=I_1\times I_2\times\ldots\times I_n$
we define the $\preceq$-convex piecewise affine mapping
$P^{(i)}:X\to Y$, letting $P^{(i)}(x)=q_{i_1}(x)e_1 +
q_{i_2}(x)e_2 + \ldots + q_{i_n}(x)e_n.$ Then
\begin{equation}
  \label{e5.6}
P(x) = \inf\limits_{i \in I}P^{(i)}(x), \ x\in X.
\end{equation}

By Proposition~\ref{pr5.2} each $P^{(i)},\ i \in I,$ can be
represented in the form (\ref{e5.2}) and, consequently, it follows
from Equality (\ref{e5.6}), $P$ can be represented in the form
(\ref{e5.3}). It proves the implication $a)\Rightarrow b).$

$b)\Rightarrow d)$ Suppose that a mapping $P:X\to Y$ is of the
form (\ref{e5.3}) and let $G_i(x)=\sup\limits_{1\le j\le
m(i)}A_{ij}(x),$ $i=1, \ldots, k.$ Then
$$
P(x) = \inf\limits_{1\le i\le k}G_i(x)=-\sup\limits_{1\le i\le k}
\left(\sum\limits_{s=1, s\ne
i}^kG_s(x)-\sum\limits_{s=1}^kG_s(x)\right)=
$$
$$
= \sum\limits_{s=1}^kG_s(x)- \sup\limits_{1\le i\le
k}\sum\limits_{s=1, s\ne i}^kG_s(x)
$$
Since a sum and a supremum of finitely many $\preceq$-convex
piecewise mappings also are $\preceq$-convex piecewise affine
mappings, both mappings $x\to \sum\limits_{s=1}^kG_s(x)$ and $x\to
\sup\limits_{1\le i \le k}\sum\limits_{s=1, s \ne i}^kG_s(x)$ are
$\preceq$-convex and piecewise affine and, consequently, each of
them can be represented in the form (\ref{e5.2}). It proves that
the mapping $P$ can be represented in the form (\ref{e5.5}).

$d)\Rightarrow c)$ Since $P$ is of the form (\ref{e5.5}),  letting
$A_{ij}(x)=C_i(x)-D_j(x),$\, $x\in X,\, j= 1, \ldots, m,\ i=1,
\ldots, k$, we immediately get the representation (\ref{e5.4}).
This completes the proof of the theorem.
\medskip

It immediately follows from the assertion d) of Theorem
\ref{th5.1} that for every piecewise affine mapping $P:X\to Y$
there exists finite two-index family of affine mappings
$\{A_{ij}:X \to Y, i=1,\ldots,m;j=1,\ldots,k\}$ such that
$$
P(x)=\inf_{1\le i\le m}\sup_{1\le j\le k}A_{ij}(x)= \sup_{1\le j\le
k} \inf_{1\le i\le m} A_{ij}(x). 
$$

For piecewise affine functions (that is, for the case when
$Y={\mathbb{R}}$) a straightforward proof of the above
representation was given in [14].

\vskip4pt plus2pt

We complete this section with the theorem that characterizes
approximative properties of piecewise affine mappings in the space
of continuous mappings.

\begin{theorem}{\bf \hspace{-6pt}.}\label{th5.2}
Let $\epsilon>0$ be an arbitrary positive real. For any continuous
mapping $F:X\to Y$ and for any compact subset $Q$ of  $X$ there
exists a piecewise affine mapping $P:X\to Y$ such that
$$
\max\limits_{x\in Q}||F(x)-P(x)||_Y < \epsilon.
$$
\end{theorem}

The proof of this theorem follows from the lattice version of the
Stone-Weirstrass theorem [29] (see, also, [30, Theorem 9.12]).

\begin{cor}{\bf \hspace{-6pt}.}\label{cor5.1}
The space $P(X, Y)$ of piecewise affine mappings is dense in the
space $C(X, Y)$ of continuous mappings from $X$ into $Y$ endowed
with the local convex topology generated by the family of seminorm
$||P||_Q=\max\limits_{x\in Q}||P(x)||$, where $Q$ runs through
compact sets in $X$.
\end{cor}
\bigskip


\section{Piecewise linear mappings}

\setcounter{equation}{0}

In this section we present the main properties of piecewise linear
mappings without proofs.

A piecewise affine mapping $P:X\to Y$ is called {\it piecewise
linear}, if it is positively homogenous, that is, if $P(\lambda x) =
\lambda P(x)$ for all $x \in X$ and $\lambda \ge 0.$

The collection of all piecewise linear mappings from $X$ into $Y$
will be denoted by $PL(X, Y)$. First of all we note that
composition of piecewise linear mappings is a piecewise linear
mapping too. If a basis is chosen and fixed in the space $Y$ then
a mapping $P:X\to Y$ is piecewise linear if and only if its
coordinate functions are piecewise linear. A mapping $P:X\to Y$ is
piecewise linear if and only if for any polyhedral partial order
$\preceq$, defined on $Y,$ $\preceq$-epigraph ${\rm epi}_\preceq
P$ and $\preceq$-hypograph ${\rm hyp}_\preceq P$ are polyhedral
cones in $X\times Y$.

If a partial order $\preceq$ is defined on $Y$ by a minihedral
convex cone or, equivalently, if $Y$ ia an Archimedean vector
lattice then with respect to standard pointwise algebraic
operations and the partial order $\preccurlyeq$ defined by
\begin{equation}
  \label{e6.1}
P\preccurlyeq Q \Longleftrightarrow P(x)\preceq Q(x) \ \mbox{ для
всех}\ x\in X,
\end{equation}
the collection $PL(X,\,Y)$ of piecewise linear mappings is a vector
sublattice of $PA(X, Y)$ (and of ${\mathcal{F}}(X,\,Y)$ as well).

A piecewise linear mapping $P:X\to Y$ is $\preceq$-convex if and
only if it can be represented in the form
$$
P(x) = \sup\limits_{j\in J}L_j(x), \ x\in X,
$$
where $L_j:X\to Y, \ j\in J,$ are linear mappings from $X$ into $Y$,
$J$ is a finite family of indices and the operation of supremum is a
lattice operation in $Y.$

\begin{theorem}{\bf \hspace{-6pt}.}\label{th7.1}
Let $P:X \to Y$ be a mapping from $X$ into $Y.$
The following assertions are equivalent:\\
a) $P:X\to Y$ is piecewise linear;\\
b) $P:X\to Y$ can be presented in the form
$$
P(x) = \inf\limits_{1\le i\le k}\sup\limits_{1\le j\le
m(i)}L_{ij}(x),\ \ x\in X,
$$
where $L_{ij}:X\to Y,\ j=1, \ldots, m(i),\ i=1, \ldots, k$, are linear mappings;\\
c) $P:X\to Y$ can be presented in the form
$$
P(x) = \sup\limits_{1\le i\le k}\inf\limits_{1\le j\le
m(i)}L_{ij}(x),\ \ x\in X,
$$
where $L_{ij}:X\to Y,\ j=1, \ldots, m(i),\ i=1, \ldots, k$ are linear mappings;\\
d) $P:X\to Y$ can be presented as difference of two
$\preceq$-convex piecewise linear mappings, that is
$$ P(x)=\sup\limits_{1\le i \le k}U_i(x) - \sup\limits_{1\le j
\le m}V_j(x),\ x\in X,
$$
where $U_i, V_j:X\to Y,\ j=1, \ldots, m,\ i=1,\ldots, k$, are linear
mappings.

In assertions b), c) and d) operations of supremum and infimum are
lattice operations in $Y.$
\end{theorem}

Let ${\cal P}(X, Y)$ be the vector space of positively homogenous
continuous mappings from $X$ into $Y.$ The space ${\cal P}(X, Y)$ is
a Banach one with respect to the norm
$$
||P||=\sup\limits_{x\in S}||P(x)||_Y,
$$
where $S$ is the unit sphere in $X.$

It follows from the lattice version of the Stone-Weierstrass
theorem [29, 30]), that the vector subspace $PL(X, Y)$ of
piecewise linear mappings is dense in ${\cal P}(X, Y)$.


\section*{Acknowledgment}

{The research was supported by the National Program of Fundamental
Researches of Belarus, the grant “Convergence — 1.4.03”.}


\bigskip
{\bf{References}}
\bigskip

\makebox[2em][r]{[1] } {C.P. Rourke,  B.J. Sanderson,}
{Introduction in Piecewise-linear Topology, Springer,  Berlin et
al, 1972.}

\makebox[2em][r]{[2] } {A.I. Subbotin,} {Minimax inequalities and
Hamilton-Jacobi equations, Nauka, Мoscow, 1991 (in Russian).}

\makebox[2em][r]{[3] } A.I. Subbotin, A piecewise-linear cost
function for a differential game with simple moves, Trudy Mat.
Inst. Steklov., 185 (1988),  242–251.

\makebox[2em][r]{[4] } {D. Melzer,} {On the expressibility of
piecewise-linear continuous functions as the difference of two
piecewise-linear convex functions, Mathematical Programming Study
29(1986)~118~--~134.}

\makebox[2em][r]{[5] } {A. Kripfganz, R. Schulze,} {Piecewise
affine functions as the difference of two convex functions,
Optimization 18(1987) 23~--~29.}

\makebox[2em][r]{[6] } {S.G. Bartels, L. Kuntz, S. Scholtes,}
{Continuous selections of linear functions and nonsmooth critical
points, Nonlinear Analysis. Theory, Methods and Applications
24(1995) 385~--~407.}

\makebox[2em][r]{[7] } {R.T. Rockafellar, R.J.-B. Wets},
Variational analysis, Springer-Verlag, Berlin, 1998.

\makebox[2em][r]{[8] } {S. Ovchinnikov,} Max-Min representation of
piecewise linear functions, Contributions to Algebra and Geometry
43(2002) 297-302.

\makebox[2em][r]{[9] } Ch.D. Aliprantis, R. Tourky, Cones and
Duality (Graduate Studies in Mathematics, V.84), Providence, RI,
American Mathematical Society, 2007.

\makebox[2em][r]{[10] } {V.V. Gorokhovik, O.I. Zorko}, {Polyhedral
quasidifferentiability of real-valued functions, Doklady Akademii
Nauk Belarusi 36(1992) 393~--~397 (in Russian).}

\makebox[2em][r]{[11] } {V.V. Gorokhovik, O.I. Zorko}, {Piecewise
affine functions and polyhedrel sets, Optimization 31(1994)
3~--~17.}

\makebox[2em][r]{[12] } {V.V. Gorokhovik, O.I. Zorko}, {Nonconvex
polyhedral sets and functions and their analytical
represen\-tations, Doklady Akademii Nauk Belarusi 39(1995) 5~--~9
(in Russian).}

\makebox[2em][r]{[13] } {V.V. Gorokhovik, O.I. Zorko}, Piecewise
affine mappings, in: N.A. Evkhuta (Ed.), Operators and operator
equations, Novocherkask State Technical University, Novocherkask,
1995, pp.~3~--~18 (in Russian).

\makebox[2em][r]{[14] } {V.V. Gorokhovik, D.S. Malashevich,} {On
analytical representations of polyhedral sets and piecewise affine
functions, The National Academy of Sciences of Belarus,
Proceedings of the Institute of Mathematics 5(2001) 15~--~19 (in
Russian).}

\makebox[2em][r]{[15] } {V.V. Gorokhovik,} {Geometrical and
analytical characterizations of piecewise affine mappings, The
National Academy of Sciences of Belarus, Proceedings of the
Institute of Mathematics 15 (2007) 22~--~32 (in Russian).}

\makebox[2em][r]{[16] } {B. Gr\"{u}nbaum,} {Convex Polytopes, John
Wiley \& Sons, London-New York-Sydney, 1967.}

\makebox[2em][r]{[17] } {R.T. Rockafellar,} {Convex analysis,
Princeton University Press, Princeton, 1970.}

\makebox[2em][r]{[18] } {K. Leichtweiss,} {Konvexe Mengen, VEB
Deitscher Verlag der Wissenschaften, Berlin, 1980.}

\makebox[2em][r]{[19] } {A. Brondsted,} {An Introduction to Convex
Polytopes, Springer-Verlag, New York-Heidelberg-Berlin, 1983.}

\makebox[2em][r]{[20] } {H.W. Kuhn, A.W. Tucker (Eds.),} {Linear
inequalities and related systems, Princeton University Press,
Princeton, 1956.}

\makebox[2em][r]{[21] } {S.N. Chernikov,} {Linear inequalities,
Nauka, Moscow, 1968 (in Russian).}

\makebox[2em][r]{[22] } {G.B. Dantzig,} {Linear programming and
extensions, Princeton University Press, Princeton, 1963.}

\makebox[2em][r]{[23] } {S.A. Ashmanov,} {Linear programming,
Nauka, Moscow, 1981(in Russian).}

\makebox[2em][r]{[24] } {S. Robinson,} {Some Continuity Properties
of Polyhedral Multifunctions, Mathematical Prog\-ram\-ming Study,
14(1981)~206~--~214.}

\makebox[2em][r]{[25] } {G. Birkhoff,} {Lattice Theory, American
Mathematical Society, Providence, R.I., 1967.}

\makebox[2em][r]{[26] } {G. Gr{a}tzer,} {General Lattice Theory,
Springer, Berlin, 1978.}

 {\Russian
\makebox[2em][r]{[27] } {N. Bourbaki} {Topologie generale. Ch. IX -- X.}-- Moscow: Nauka, 1975
(Russian traslation). Использование вещественных чисел в общей топологии. Функциональные
пространства. Сводка результатов. Словарь.}

\makebox[2em][r]{[28] }  I.M. Glazman, Yu.I. Lyubich,
Finite-dimensional linear analysis: a systematic presentation in
problem form, M.I.T. Press, Cambridge, Mass. 1974.

\makebox[2em][r]{[29] } {N. Bourbaki} {Integration. Measures,
integration of measures}

\makebox[2em][r]{[30] } {Ch.D. Aliprantis, K.C. Border,} {Infinite
dimensional analysis, Berlin, Springer, 1999.}

\end{document}